\documentclass{article}
\usepackage{mystyle}
\begin{document}
\title{Usuba's extendible cardinal}
\author{Gabriel Goldberg}
\maketitle
\begin{abstract}
    Answering a question of Usuba \cite{UsubaExtendible}, we show that an extendible cardinal 
    can be preserved by a set forcing that is not a small forcing. 
\end{abstract}
\section{Introduction}
This note addresses the following question of Usuba \cite{UsubaExtendible} by showing
that its answer is independent of ZFC: 
\begin{qst}\label{qst:usuba}
    Suppose \(W\) is a set forcing extension of \(V\)
    and \(\kappa\) is a cardinal that is extendible in \(W\).
    Must there exist a partial order \(\mathbb P\in V\)
    of cardinality less than \(\kappa\) and
    a \(V\)-generic filter \(G\) on \(\mathbb P\)
    such that \(W = V[G]\)?
\end{qst}

Usuba's question is motivated by his celebrated theorem
on the relationship between extendible cardinals and forcing.
Recall that if \(M\subseteq N\) are models of ZFC,
then \(M\) is a \textit{ground of \(N\)} if
\(N\) is a set forcing extension of \(M\).
Surprisingly, the grounds of a model of ZFC can be defined uniformly
over that model, enabling one to study
the structure of the grounds of the 
universe of set theory itself, which are usually
simply called \textit{grounds}.
\begin{thm*}[{Usuba \cite{UsubaExtendible}}]
    If there is an extendible cardinal,
    the universe of sets has a minimum ground.
\end{thm*}
Another theorem of Usuba shows that 
the intersection of all grounds of \(V\) is in general a model
of ZFC, which is called the \textit{mantle}
and denoted by \(\mathbb M\). If there is a minimum ground,
then it is of course equal to the mantle. 
At heart, \cref{qst:usuba} asks about 
the distance between the universe of sets and the mantle:
\begin{qst*}
    Suppose \(\kappa\) is an extendible cardinal.
    Must there be a partial order \(\mathbb P \in V_\kappa\cap \mathbb M\)
    and an \(\mathbb M\)-generic filter \(G\subseteq \mathbb P\)
    such that \(V = \mathbb M[G]\)?
\end{qst*}
If no such partial order exists, then Usuba's Theorem fails in \(V_\kappa\)
where \(\kappa\) is the least extendible cardinal. Since one would expect \(V_\kappa\)
to satisfy every large cardinal axiom short of extendibility,
the consistency of a negative
answer to this question would argue that the large cardinal
hypothesis of Usuba's Theorem cannot be weakened.

The \textit{Ground Axiom}, introduced by Hamkins-Reitz, states that \(V = \mathbb M\). 
A consequence of Usuba's Theorem is that if there is an extendible cardinal,
the mantle satisfies the Ground Axiom.
The \textit{Local Ground Axiom}
states that the Ground Axiom holds in 
\(V_\lambda\) whenever \(\lambda\) is a Beth fixed point.
It is not hard to force the Local Ground Axiom by a 
class Easton iteration that preserves extendible cardinals.
\begin{repprp}{prp:positive}
    Assume that the mantle satisfies the Local Ground Axiom.
    Then the answer to \cref{qst:usuba} is yes.
\end{repprp}

More interestingly, it is also consistent that the answer to \cref{qst:usuba} is no.
A cardinal \(\lambda\) is \textit{indestructibly \(\Sigma_n\)-correct}
if \(V_\lambda\) is a \(\Sigma_n\)-elementary substructure of
any \(\lambda\)-directed closed forcing extension of \(V\).
Every \(\Sigma_1\)-correct cardinal is indestructibly \(\Sigma_1\)-correct;
every Laver indestructible supercompact cardinal is indestructibly \(\Sigma_2\)-correct;
and by a theorem of Bagaria-Hamkins-Tsaprounis-Usuba \cite{BHTU},
there cannot be an indestructibly \(\Sigma_3\)-correct cardinal.
\begin{repthm}{thm:main}
    Assume there is a proper class
    of indestructibly \(\Sigma_2\)-correct cardinals.
    Then the answer to \cref{qst:usuba} is no.
\end{repthm}

Given a model with an extendible cardinal and a proper class of supercompact cardinals,
one can force to obtain a model satisfying the hypothesis of \cref{thm:main}
using the global Laver preparation, as described in \cite[Section 1]{ApterLaverIndestructibility}. 
\cref{prp:preparation} shows that if one just wants indestructible
\(\Sigma_2\)-correctness, the supercompact cardinals are overkill: 
it suffices to start with a model
with an extendible cardinal and a proper class
of inaccessible \(\Sigma_2\)-correct cardinals.
It seems likely that the inaccessible \(\Sigma_2\)-correct cardinals
can be dispensed with as well,
but the author does not see how to do this.
\section{The Local Ground Axiom}
\subsection{Extendibility over the mantle}
In this section, we show that it is consistent
that the answer to \cref{qst:usuba} is yes.
This is straightforward and follows Usuba's original, but the proof suggests
the route to proving the much more interesting negative consistency
result by imposing two fundamental constraints on the structure of a counterexample.

If \(M\subseteq N\) are models of ZFC,
then \(M\) is a \textit{\(\kappa\)-ground of \(N\)} if there is a partial order
\(\mathbb P\in M\) of \(M\)-cardinality less than \(\kappa\)
carrying a \(M\)-generic filter \(G\) such that
\(N = M[G]\). The \textit{\(\kappa\)-mantle} is the intersection of all \(\kappa\)-grounds of \(V\).
We will need Usuba's Theorem in the following form:
\begin{thm}[Usuba]\label{thm:usuba}
    If \(\kappa\) is extendible, then the mantle is equal to the \(\kappa\)-mantle.\qed
\end{thm}

If \(M\) is an inner model of ZFC,
a cardinal \(\kappa\) is \textit{extendible over \(M\)} if
for all ordinals \(\lambda \geq \kappa\), for some \(\lambda' > \lambda\),
there is an elementary embedding \(j : V_{\lambda+1}\to V_{\lambda'+1}\)
such that \(\crit(j) = \kappa\), \(j(\kappa) > \lambda\),
\(j\restriction V_\lambda^M\) belongs to
\(M\), and \(j(V_\lambda^M) = V_{\lambda'}^M\). 
The following proof is based on Usuba's original
proof of his theorem using hyperhuge cardinals \cite{Usuba}.
\begin{prp}\label{prp:extendible_over}
    An extendible cardinal \(\kappa\) is extendible over a ground \(W\)
    if and only if \(W\) is a \(\kappa\)-ground.
    \begin{proof}
        Let \(\mathbb Q\in W\) be a partial order
        carrying a \(W\)-generic filter \(H\) such that \(V = W[H]\).

        Suppose \(\lambda\) is a Beth fixed point larger than the rank of \(\mathbb Q\).
        Let \(j : V_{\lambda+1}\to V_{\lambda'+1}\) be an elementary embedding
        such that \(\crit(j) = \kappa\), \(j(\kappa) > \lambda\),
        \(j\restriction W_\lambda\) belongs to
        \(W\), and \(j(V_\lambda^W) = V_{\lambda'}^W\).
        Since \(V_{\lambda'}^W[H] = V_{\lambda'}\),
        \(V_{\lambda'}\) satisfies that there is a partial order
        \(\mathbb P\in W_{j(\kappa)}\) carrying a \(V_{\lambda'}^W\)-generic
        filter \(G\) such that \(V_{\lambda'}^W[G] = V_{\lambda'}\).
        Therefore \(V_\lambda\) satisfies that there is a partial order
        \(\mathbb P\in V_{\kappa}^W\) carrying a \(V_{\lambda}^W\)-generic
        filter \(G\) such that \(V_{\lambda}^W[G] = V_{\lambda}\).

        By the pigeonhole principle, there is a partial order
        \(\mathbb P\in V_\kappa^W\) carrying a filter \(G\)
        such that for a proper class of Beth fixed points \(\lambda\),
        \(G\) is \(V_\lambda^W\)-generic and \(V_\lambda = V_\lambda^W[G]\).
        It follows that \(G\) is \(W\)-generic and \(V = W[G]\),
        which proves the proposition.
    \end{proof}
\end{prp}
\begin{cor}\label{cor:extendible_over_mantle}
    If \(\kappa\) is the least extendible cardinal, the mantle is a \(\kappa\)-ground
    if and only if \(\kappa\) is extendible over \(\mathbb M\).\qed
\end{cor}
The \textit{Local Ground Axiom}
states that the Ground Axiom holds in 
\(V_\lambda\) whenever \(\lambda\) is a Beth fixed point.
\begin{prp}\label{prp:lga}
    If the mantle satisfies the Local Ground Axiom and \(\kappa\) is an
    extendible cardinal, then \(\kappa\) is
    extendible over the mantle.
    \begin{proof}
        We claim that for all Beth fixed points \(\lambda > \kappa\),
        \(V_\lambda^\mathbb M = \mathbb M^{V_\lambda}\). By Usuba's Theorem \cite{UsubaExtendible},
        \(\mathbb M^{V_\lambda}\subseteq \mathbb M\). Conversely, it \(W\) is a ground of \(V_\lambda\),
        then since \(\mathbb M\) and \(W\) are grounds of \(V_\lambda\),
        they have a common ground \(N\), and since \(V_\lambda^\mathbb M\)
        satisfies the Ground Axiom, \(N = V_\lambda^\mathbb M\),
        and so \(V_\lambda^\mathbb M\subseteq W\).

        Now suppose \(\lambda\) is a Beth fixed point, and let
        \(j : V_{\lambda+1}\to V_{\lambda'+1}\) be an elementary embedding
        such that \(\crit(j) =\kappa\) and \(j(\kappa) > \lambda\).
        We must show that \(j\restriction V_\lambda^\mathbb M\) belongs to \(\mathbb M\)
        and \(j(V_\lambda^\mathbb M) = V_{\lambda'}^\mathbb M\). 
        That \(j\restriction V_\lambda^\mathbb M\) belongs to \(\mathbb M\) is
        immediate by Usuba's Theorem and the absoluteness of extendibility
        to \(\kappa\)-grounds. On the other hand, clearly \(\lambda'\) is a Beth fixed point,
        and so 
        \[j(V_\lambda^\mathbb M) = j(\mathbb M^{V_\lambda}) = \mathbb M^{V_{\lambda'}} = V_{\lambda'}^\mathbb M\qedhere\]
    \end{proof}
\end{prp}

\begin{prp}\label{prp:positive}
    If the mantle satisfies the Local Ground Axiom, then the answer to \cref{qst:usuba}
    is yes.
    \begin{proof}
        Suppose \(N\) is a set forcing extension and \(\kappa\) is
        extendible in \(N\). Then applying \cref{prp:lga} in \(N\),
        \(\kappa\) is extendible over the mantle, and so by \cref{cor:extendible_over_mantle},
        the mantle is a \(\kappa\)-ground of \(N\). 
        Since \(\mathbb M\subseteq V\subseteq N\), the Intermediate Model Theorem implies
        that \(V\) is a \(\kappa\)-ground of \(N\).
    \end{proof}
\end{prp}
\subsection{Two constraints}
\cref{prp:positive} highlights two key constraints
that guide the way to a counterexample to Usuba's question.
The first constraint is that by \cref{prp:positive}, it is consistent that
\cref{qst:usuba} has a positive answer in all forcing extensions.
This suggests that in order to find a counterexample,
one should start with a preparatory class forcing.

The second constraint is more subtle.
By \cref{cor:extendible_over_mantle},
if one is to preserve an extendible cardinal
by a set forcing that is not a small forcing, this preservation cannot be proved
use the standard lifting arguments for extendible cardinals.
The reason is that if a cardinal \(\kappa\)
is shown to be extendible in a forcing extension \(N\) of a model \(M\) using these lifting arguments, 
then \(N\) satisfies that \(\kappa\) is extendible over \(M\), and hence \(M\) is a \(\kappa\)-ground of \(N\).

Since it is hard to see how to preserve extendible cardinals
without a lifting argument, answering \cref{qst:usuba} seems to
require inventing a novel forcing notion
along with an entirely new preservation argument for extendible cardinals, a
task that lies above the author's paygrade.
The solution instead is simply to
\textit{reformulate} extendibility in terms of normal fine ultrafilters (\cref{lma:sigma2_char})
and then to employ the standard lifting arguments
from the theory of supercompactness to show that, thanks
to our preliminary preparatory forcing, this reformulation
is preserved by a forcing notion that is about as far from 
novel as one can get: an Easton product of Cohen forcings.
\section{Indestructible correctness}
\subsection{The main theorem}
In this section, we show that under indestructibility
hypotheses, the answer to \cref{qst:usuba} is no.
A cardinal \(\lambda\) is \textit{\(\Sigma_n\)-correct}
if \(V_\lambda\) is a \(\Sigma_n\)-elementary substructure of \(V\).
A cardinal \(\lambda\) is \textit{indestructibly \(\Sigma_n\)-correct}
if \(V_\lambda\) is a \(\Sigma_n\)-elementary substructure of \(V\)
in any \(\lambda\)-directed closed forcing extension.
\begin{thm}\label{thm:main}
    Assume there is a proper class of indestructibly \(\Sigma_2\)-correct
    cardinals. Then the answer to \cref{qst:usuba} is no.
\end{thm}
\cref{prp:preparation} shows that the hypothesis of \cref{thm:main}
can be class forced starting with a model with a proper class
of inaccessible \(\Sigma_2\)-correct cardinals. 

The key to the proof of \cref{thm:main} 
is the following characterization of extendibility,
observed independently and earlier by Bagaria.
For any cardinal \(\lambda\), let \(T_{\kappa,\lambda}\) (resp.\ \(T_{\kappa,\lambda}^*\))
denote the set of \(\sigma\in P_\kappa(\lambda)\) such that 
the ordertype of \(\sigma\) is \(\Sigma_2\)-correct (resp.\ indestructibly
\(\Sigma_2\)-correct) in \(V_\kappa\).
In most cases of interest, \(\kappa\) itself will be \(\Sigma_2\)-correct,
in which case a cardinal is \(\Sigma_2\)-correct in \(V_\kappa\) if and only if
it is truly \(\Sigma_2\)-correct. Similarly, if \(\kappa\) is \(\Sigma_2\)-correct,
then a cardinal is indestructibly \(\Sigma_2\)-correct in \(V_\kappa\) if and only if
it is truly indestructibly \(\Sigma_2\)-correct.
\begin{lma}\label{lma:sigma2_char}
    A cardinal \(\kappa\) is extendible if and only if for arbitrarily large
    cardinals \(\lambda\), there is a normal fine \(\kappa\)-complete ultrafilter
    on \(T_{\kappa,\lambda}\).
    \begin{proof}
        For the forwards direction, suppose \(\lambda\) is \(\Sigma_2\)-correct
        and \(\kappa\) is \(\lambda\)-extendible, and we will show that there is a 
        normal fine \(\kappa\)-complete ultrafilter
        on \(T_{\kappa,\lambda}\). Let \(j : V_{\lambda+1}\to V_{\lambda'+1}\) be an elementary embedding
        with critical point \(\kappa\). Since \(\lambda\) is \(\Sigma_2\)-correct,
        \(\lambda\) is \(\Sigma_2\)-correct in \(V_{j(\kappa)}\),
        and so \(j[\lambda]\in T_{j(\kappa),j(\lambda)}\). It follows that
        there is a normal fine \(\kappa\)-complete ultrafilter on 
        \(T_{\kappa,\lambda}\); namely, the ultrafilter derived from \(j\) using \(j[\lambda]\),
        or in symbols,
        \(\{A\subseteq T_{\kappa,\lambda} : j[\lambda]\in j(A)\}\).

        Now we show that if there is a normal fine \(\kappa\)-complete ultrafilter \(\mathcal U\)
        on \(T_{\kappa,\lambda}\), then \(\kappa\) is \(\gamma\)-extendible for all \(\gamma < \lambda\).
        In particular, this implies the reverse direction of the lemma.
        To see this, let \(j : V\to M\) be the ultrapower embedding associated to \(\mathcal U\).
        Note that \(j\restriction V_\lambda\) belongs to \(V_{j(\kappa)}\cap M\), and in \(M\),
        it is an elementary embedding from \(V_\lambda\) to \(V_{\lambda'}^M\). 
        Therefore \(V_{j(\kappa)}^M\) satisfies that \(\kappa\) is \(\gamma\)-extendible
        for all \(\gamma < \lambda\). 
        Since \(V_\lambda\) is a \(\Sigma_2\)-elementary substructure of \(V_{j(\kappa)}^M\),
        \(V_\lambda\) satisfies that \(\kappa\) is \(\gamma\)-extendible for all \(\gamma < \lambda\),
        and since the \(\gamma\)-extendibility of \(\kappa\) is expressed by a \(\Sigma_2\)-formula, 
        this is upwards absolute to \(V\). Hence \(\kappa\) is 
        is \(\gamma\)-extendible for all \(\gamma < \lambda\), as claimed.
    \end{proof}
\end{lma}

The proof of \cref{lma:sigma2_char} yields:
\begin{lma}\label{lma:indestructible_sigma2_char}
    If \(\kappa\) is extendible and \(\lambda\) is indestructibly \(\Sigma_2\)-correct,
    then there is a normal fine \(\kappa\)-complete ultrafilter
    on \(T_{\kappa,\lambda}^*\).\qed
\end{lma}
Given this, we turn to the proof of \cref{thm:main}
\begin{proof}[Proof of \cref{thm:main}]
        For each cardinal \(\delta\),
        let \(\lambda_\delta\) denote the least indestructibly \(\Sigma_2\)-correct cardinal \(\lambda\)
        such that there is no normal fine \(\delta\)-complete ultrafilter on \(T_{\kappa,\lambda}^*\)
        if there is such a cardinal.
        Let \(\mathbb P_\delta = \text{Add}(\lambda_\delta^+,1)\),
        and let \(\mathbb P\) be the Easton product of the partial orders \(\mathbb P_\delta\)
        for \(\delta < \kappa\). 

        Let \(G\subseteq \mathbb P\) be a \(V\)-generic filter. Clearly there is no partial order
        \(\bar{\mathbb P}\in V_\kappa\) carrying a \(V\)-generic filter
        \(\bar G\in V[G]\) such that \(V[G] = V[\bar G]\).
        To show that the answer to \cref{qst:usuba} is no, it therefore suffices to show that
        \(\kappa\) is extendible in \(V[G]\). To prove this, we will verify
        the criterion of \cref{lma:sigma2_char}. 
        
        Note that the class of indestructibly \(\Sigma_2\)-correct cardinals is closed.
        Therefore there is a proper class of indestructibly \(\Sigma_2\)-correct \textit{singular} cardinals
        \(\lambda\), which have the property that \(2^\lambda = \lambda^+\)
        by Solovay's Theorem \cite{Solovay} that the Singular Cardinals Hypothesis holds above a supercompact cardinal.
        Suppose \(\lambda > \kappa\) is
        an indestructibly \(\Sigma_2\)-correct cardinal such that \(2^\lambda = \lambda^+\). 
        We claim that \(V[G]\) satisfies that there is a
        normal fine \(\kappa\)-complete ultrafilter on \((T_{\kappa,\lambda})^{V[G]}\).

        Working in \(V\), let \(\mathcal U\) be Mitchell minimal
        among all \(\kappa\)-complete normal fine ultrafilters 
        on \(T_{\kappa,\lambda}^*\). Such an ultrafilter exists by \cref{lma:indestructible_sigma2_char}
        and the wellfoundedness of the Mitchell order.
        We claim that \(\lambda = (\lambda_\kappa)^{M_\mathcal U}\).
        To see this, note that since \(V_\kappa^{M_\mathcal U} = V_\kappa\),
        for all \(\gamma \leq \lambda\), \((T_{\kappa,\gamma}^*)^{M_\mathcal U} = T_{\kappa,\gamma}^*\).
        Moreover, since \(V_\lambda\subseteq M_{\mathcal U}\) and \(\lambda\) is a strong limit cardinal,
        every ultrafilter on \(T_{\kappa,\gamma}^*\) for \(\gamma < \lambda\) belongs to \(M_\mathcal U\).
        Therefore \((\lambda_\kappa)^{M_\mathcal U} \leq \lambda\). But
        there can be no normal fine ultrafilter \(\mathcal W\) on \(T_{\kappa,\lambda}^*\) in \(M_\mathcal U\):
        otherwise \(\mathcal W\) is a normal fine ultrafilter on \(T_{\kappa,\lambda}^*\) in \(V\)
        since \(P(T_{\kappa,\lambda})\subseteq M_\mathcal U\), and this contradicts the Mitchell minimality
        of \(\mathcal U\). It follows that \(\lambda \geq (\lambda_\kappa)^{M_\mathcal U}\), which proves the claim.

        Since \(\lambda = (\lambda_\kappa)^{M_\mathcal U}\), the forcing \(j_\mathcal U(\mathbb P)\)
        is isomorphic to the product \(\mathbb P\times \mathbb Q\) where 
        \(\mathbb Q = \prod_{\kappa \leq \delta < j_\mathcal U(\kappa)} \mathbb P_\delta\)
        is \(\lambda^+\)-directed closed in \(M_\mathcal U\). Since \(M_\mathcal U\) is closed under
        \(\lambda\)-sequences, \(\mathbb Q\) really is \(\lambda^+\)-directed closed. 
        In addition, \(|P^M(\mathbb Q)|^M \leq (2^{\kappa})^{M_\mathcal U} < j_\mathcal U(\lambda) < (2^\lambda)^+ = \lambda^{++}\).
        The final bound follows from the fact that \(2^\lambda = \lambda^+\).
        Therefore in \(V\), \(|P^M(\mathbb Q)| \leq \lambda^+\) 
        and \(\mathbb Q\) is \(\lambda^+\)-closed, and so one can build an \(M\)-generic
        filter \(H\subseteq \mathbb Q\) with \(H\in V\). 

        The closure of \(\mathbb Q\) implies that 
        \(M[H]\) contains no new dense subsets of \(\mathbb P\), and so \(G\) is an
        \(M[H]\)-generic filter on \(\mathbb P\).
        By standard results on mutual genericity,
        this means that \(G\times H\) is an \(M\)-generic filter on \(\mathbb P\times \mathbb Q\).
        The cardinal \(\lambda\) is \(\Sigma_2\)-correct in \(M[H]\) by since
        \(\lambda\) is indestructibly \(\Sigma_2\)-correct in \(M\).
        Since \(G\) is \(M[H]\)-generic for a forcing in \((V_\lambda)^{M[H]}\),
        \(\lambda\) is \(\Sigma_2\)-correct in \(M[H\times G]\).

        Finally, identifying \(j_\mathcal U(\mathbb P)\) with \(\mathbb P\times \mathbb Q\) in the natural way,
        \(j_\mathcal U[G] = G\times \{1\}\subseteq G\times H\), and so the embedding \(j_\mathcal U\)
        lifts uniquely to an elementary embedding \(j : V[G]\to M[G\times H]\) such that
        \(j(G) = G\times H\).
        Since \(\lambda\) is \(\Sigma_2\)-correct in \(M[G\times H]\), the set
        \(j[\lambda]\) belongs to \(T_{j(\kappa),j(\lambda)}\) as computed in \(M[G\times H]\).
        As a consequence, working in \(V[G]\), 
        the ultrafilter 
        \[\mathcal D = \{A\subseteq T_{\kappa,\lambda} : j[\lambda]\in j(A)\}\]
        derived from \(j\) using \(j[\lambda]\) 
        is a normal fine \(\kappa\)-complete ultrafilter on \(T_{\kappa,\lambda}\). 
\end{proof}
\subsection{Forcing indestructible correctness}
In this section, we show that the
hypotheses of \cref{thm:main} are consistent relative to a large
cardinal hypothesis just past an extendible cardinal. 
\begin{prp}\label{prp:preparation}
    There is a class forcing \(\mathbb Q\) such that
    every inaccessible \(\Sigma_2\)-correct cardinal of \(V\)
    is an inaccessible indestructibly \(\Sigma_2\)-correct cardinal of \(V^\mathbb Q\).
    Moreover this forcing preserves extendible cardinals.
\end{prp}
\begin{proof}
Suppose \(\eta\) is an ordinal, \(y\) is a set of rank at most \(\eta\),
and \(\varphi(x)\) is a set-theoretic formula in one free variable. Let \(F(\eta,\varphi,y)\) 
denote the least
Beth fixed point \(\beta \geq \rank(y)\) for which there is an
\(\eta\)-directed closed partial order
\(\mathbb P\in V_\beta\) such that \(V_\beta^\mathbb P\vDash \varphi(y)\).
For each ordinal \(\eta\), let
\(f(\eta)\ = \sup_{\varphi,y\in V_\eta} F(\eta,\varphi,y)\). 
Let \(\mathbb P(\eta)\) denote the
lottery sum of all \(\eta\)-directed closed forcings in \(V_{f(\eta)}\).

We define a class Easton iteration 
\(\langle \mathbb Q_\alpha,\dot{\mathbb P}_\alpha : \alpha \in \Ord\rangle\)
and a continuous sequence of ordinals \(\eta_\alpha\)
by letting \(\dot {\mathbb P}_\alpha\) be the canonical name for
\(\mathbb P(\eta_\alpha)\) as computed in \(V^{\mathbb Q_\alpha}\)
and letting 
\(\eta_{\alpha+1}\) be the least Beth fixed point
above \(\rank(\mathbb P(\eta_\alpha))\).
Let \(\mathbb Q\) be the class direct limit of the iteration and let \(G\)
be a \(V\)-generic filter on \(\mathbb Q\). 
If \(\kappa\) is an inaccessible \(\Sigma_2\)-correct cardinal,
then \(\kappa\) is closed under \(f\), and so \(\mathbb Q^{V_\kappa} \cong \mathbb Q_\kappa\).
In particular, \(\mathbb Q_\kappa\) is \(\kappa\)-cc, and so
\(\kappa\) remains inaccessible in \(V[G]\).

We will show that every \(\Sigma_2\)-correct inaccessible cardinal \(\kappa\) of \(V\) becomes 
indestructibly \(\Sigma_2\)-correct in \(V[G]\).
Assume that in \(V[G]\),
\(\beta > \kappa\) is an ordinal and there is a \(\kappa\)-directed closed partial
order \(\mathbb P\in V[G]_\beta\) such that 
\(V[G]_\beta^{\mathbb P}\vDash \varphi(y)\) where \(\varphi\) is a formula and
\(y\in V[G]_\kappa\). 
We claim that there is some \(\bar \beta < \kappa\) such that
\(V[G]_{\bar \beta}\vDash \varphi(y)\).
For ease of notation, and without loss of generality, we assume \(y\in V_\kappa\).
(More formally, we work in \(V[G\restriction \lambda]\) where \(\lambda < \kappa\)
is the hereditary cardinality of \(y\).)

For any ordinal \(\xi\), let \(D_\xi\) denote the class of conditions \(p\in \mathbb Q\) that force 
\(V^\mathbb Q_\xi\vDash \varphi(y)\), and let
\(C_\xi\) be the class of conditions \(p\in \mathbb Q\) such that 
for some \(\dot{\mathbb P}\in V^\mathbb Q_\xi\), \(p\) forces 
\((V^\mathbb Q)^{\dot{\mathbb P}}_\xi\vDash \varphi(y)\). 
We will show that for any \(p\in C_\beta\),
for some \(\bar \beta < \kappa\),
there is a condition \(r\in D_{\bar \beta}\) compatible with \(p\).
It follows that \(\bigcup_{\bar \beta < \kappa} D_{\bar \beta}\) is predense below any element of \(C_\beta\).
Since \(V[G]_\beta^{\mathbb P}\vDash \varphi(y)\),
there is a condition \(q\in G\cap C_\beta\), and so since \(\bigcup_{\bar \beta < \kappa} D_{\bar \beta}\) is
predense below \(q\), 
\(G\cap \bigcup_{\bar \beta < \kappa} D_{\bar \beta}\) is nonempty.
In other words, there is some \(\bar \beta < \kappa\) such that
\(V[G]_{\bar \beta}\vDash \varphi(y)\).

Fix \(p\in C_\beta\).
Let \(\bar p = p\restriction \kappa\in \mathbb Q_\kappa\).
Let \(\gamma = \text{supp}(\bar p)\) be the support of \(\bar p\).
Note that if \(H\) is a \(V\)-generic filter on \(\mathbb Q\)
containing \(p\), then 
in \(V[H\restriction \gamma]\), one can find an ordinal \(\bar \beta\) and an
\(\eta_\gamma\)-closed partial order \(\mathbb P\in V[H\restriction \gamma]_{\bar \beta}\) such that
\(V[H\restriction \gamma]^{\mathbb P}_{\bar \beta}\vDash \varphi(y)\); namely, \(\bar \beta = \beta\)
and \(\mathbb P = \mathbb Q_{\gamma,\beta}*\dot{\mathbb P}\) where
\(\dot{\mathbb P}\) is a \(\mathbb Q_{\gamma,\beta}\)-name for a \(\kappa\)-directed closed
partial order \(\mathbb P\) in \(V[H\restriction \beta]\) such that \(V[H]_\beta^{\mathbb P}\vDash \varphi(y)\).
The least such ordinal \(\bar\beta\) lies below \(f(\eta_\gamma)\) as computed in
\(V[H\restriction \gamma]\), and so since \(\kappa\) is \(\Sigma_2\)-correct
in \(V[H\restriction \gamma]\), \(\bar \beta < \kappa\).

It follows from the definition of
\(\dot{\mathbb P}_\gamma\) that there is an extension \(r\leq \bar p\) in \(\mathbb Q_{\gamma+1}\)
forcing that \(V^{\mathbb Q_{\gamma+1}}_{\bar \beta}\vDash \varphi(y)\).
By construction, \(V^{\mathbb Q_{\gamma+1}}_{\bar \beta} = V^{\mathbb Q}_{\bar \beta}\),
and so \(r\) forces that \(V^{\mathbb Q}_{\bar \beta}\vDash \varphi(y)\), or
in other words, \(r\in D_{\bar \beta}\).
Note that \(r\) and \(p\) are compatible since \(r\) extends \(p\restriction \kappa\)
and \(\text{supp}(r) < \kappa\). 

Finally, we sketch a proof that the forcing \(\mathbb Q\) preserves extendible cardinals.
Suppose \(\kappa\) is extendible and \(\lambda > \kappa\) is a 
limit of \(\Sigma_2\)-correct inaccessible cardinals. Also assume that
\(2^\lambda = \lambda^+\). (By Solovay's Theorem \cite{Solovay}, there is
a proper class of cardinals \(\lambda\) with these properties.)
Let \(\mathcal U\) be a normal fine \(\kappa\)-complete ultrafilter
on \(T_{\kappa,\lambda}\). Using a master condition and the fact that \(2^\lambda = \lambda^+\),
in \(V[G]\), \(j_\mathcal U : V\to M_\mathcal U\) lifts to an elementary embedding
\(j : V[G]\to M_\mathcal U[H]\). We have shown that
every inaccessible \(\Sigma_2\)-correct cardinal is \(\Sigma_2\)-correct in
\(V[G]\), and so by elementarity, \(\lambda\) remains a limit of \(\Sigma_2\)-correct cardinals
in \(M[H]\). It follows that \(\lambda\) is 
\(\Sigma_2\)-correct in \(M[H]\),
and therefore one can derive from \(j\) a normal fine \(\kappa\)-complete \(V[G]\)-ultrafilter
on \((T_{\kappa,\lambda})^{V[G]}\). 
\end{proof}
\bibliographystyle{plain}
\bibliography{Bibliography}
\end{document}